\documentclass[a4paper]{article}  
\usepackage{amsmath, amstext, amsxtra, amsthm, amscd, amsgen,
 amsbsy, amsopn, amsfonts, latexsym, amssymb, amscd}

\begin{document}



\def\proofend{\hbox to 1em{\hss}\hfill $\blacksquare
$\bigskip } 
\def\powser#1{\lbrack \lbrack #1 \rbrack \rbrack }  
  
\newtheorem{theorem}{Theorem}[section]  
\newtheorem{proposition}[theorem]{Proposition}  
\newtheorem{lemma}[theorem]{Lemma}  
\newtheorem{remark}[theorem]{Remark}  
\newtheorem{remarks}[theorem]{Remarks}  
\newtheorem{definition}[theorem]{Definition}  
\newtheorem{corollary}[theorem]{Corollary}  
\newtheorem{example}[theorem]{Example}  
\newtheorem{assumption}[theorem]{Assumption}  
\newtheorem{problem}[theorem]{Problem}  
\newtheorem{question}[theorem]{Question}  
\newtheorem{conjecture}[theorem]{Conjecture}
\newtheorem{rigiditytheorem}[theorem]{Rigidity Theorem}  
\newtheorem{mainlemma}[theorem]{Main Lemma}
\newtheorem{claim}[theorem]{Claim}
 
\def\hra{\hookrightarrow}

\def\almostisometric{finite-order-isometric}
\def\rfpd{restricted fixed point dimension}

\def\Z{{\mathbb Z}}  
\def\R{{\mathbb R}}  
\def\Q{{\mathbb Q}}  
\def\C{{\mathbb C}}  
\def\N{{\mathbb N}}  
\def\H{{\mathbb H}}  
\def\Zp #1{{\mathbb Z }/#1{\mathbb Z}}  
\def\cpt{compact}  
\def\wt{wt}
\def\cowt{co\wt} 

\def\codim{{\rm{codim}\ }}
  
\def\b{bun\-dle}  
\def\pb{principal \b }  
\def\vb{vector \b }  
  
\def\mfd{manifold}
\def\LFF{Lefschetz fixed point formula}  
\def\isorank{symrank}  

\def\ell{\Phi}
\def\eell{\widetilde \ell }  
\def\ddelta {{\widetilde \delta }}  
\def\eepsilon{{\widetilde \epsilon }}
\def\CC{C_0}
\def\oorder{o}   
\def\oha{\cal H}  
  
\def\paperref#1#2#3#4#5#6{\text{#1:} #2, {\em #3} {\bf#4} (#5)#6}  
\def\bookref#1#2#3#4#5#6{\text{#1:} {\em #2}, #3 #4 #5#6}  
\def\preprintref#1#2#3#4{\text{#1:} #2 #3 (#4)}

\hyphenation{man-i-fold equiv-a-riant in-te-ger mod-u-lo tor-sion
re-pre-sen-ta-tion di-men-sion-nal}


\title{Obstructions to positive curvature and symmetry}  
\author{Anand Dessai}
\date{}
\maketitle  
  
\begin{abstract}  
  
\noindent
We show that the indices of certain twisted Dirac operators vanish on a $Spin$-manifold $M$ of 
positive sectional curvature if the symmetry rank of $M$ is $\geq 2$ or if the symmetry rank is one 
and $M$ is two connected. We also give examples of simply connected manifolds of positive Ricci 
curvature which do not admit a metric of positive sectional curvature and positive symmetry rank.\footnote{This paper replaces the preprint ``On the elliptic genus of positively curved manifolds with symmetry'' (April 2001) which appeared at the arXiv (http://arxiv.org/abs/math.DG/0104256).}
\end{abstract}
  
\noindent
\section{Introduction}\label{intro}
An important application of index theory in Riemannian geometry is in the study of manifolds of 
positive scalar curvature. Soon after Atiyah and Singer proved the index theorem, Lichnerowicz used a 
Bochner type formula to show that the index of the Dirac operator vanishes on closed $Spin$-manifolds 
of positive scalar curvature. 

Whereas the relation between index theory and positive scalar curvature (for high dimensional simply 
connected manifolds) is well understood \cite{GrLa,ShYau1,St1,Miy} possible relations to stronger 
curvature conditions such as positive Ricci or positive sectional curvature remain obscure (see 
however the fascinating conjecture in \cite{St2}). 

In this paper we give obstructions to metrics of positive sectional curvature (positive curvature for 
short) with symmetry. We show that the indices of certain twisted Dirac operators vanish on a 
positively curved closed $Spin$-manifold $M$ provided that the symmetry rank (i.e. the rank of the 
isometry group of $M$) is at least two and the dimension of $M$ is sufficiently large. These indices 
occur as coefficients in an expansion of the elliptic genus. A similar result holds if $M$ is $2$-connected and the symmetry rank is $\geq 1$.

The elliptic genus $\varphi $ is a ring homomorphism from the oriented bordism ring to the ring of 
modular forms for 
$$\Gamma _0(2):=\{A\in SL_2(\Z )\; \mid \; A\equiv (\begin{smallmatrix} *& *\\0 &
*\end{smallmatrix})\bmod 2\}$$
(in particular, $\varphi $ vanishes in all dimensions not divisible by $4$). On the complex 
projective spaces $\C P^{2k}$ it is given by 
$$\sum _{k\geq 0}\varphi (\C P^{2k})t^{2k}=(1-2\delta \cdot t^2+\epsilon \cdot t^4)^{-1/2},$$
where $\delta $ and $\epsilon $ are modular forms of weight $2$ and $4$, respectively. 

The normalized elliptic genus $\Phi (M):=\varphi (M)/\epsilon 
^{k/2}$ of an oriented $4k$-dimen\-sional manifold $M$ expands in one of the cusps of $\Gamma _0(2)$ as a series of 
twisted signatures. Following Witten \cite{Wi} this series is best thought of as the index of a 
hypothetical signature operator on the free loop space of $M$. In the other cusp of $\Gamma _0(2)$ 
the elliptic genus expands as a series $\Phi 
_0(M)$ of characteristic numbers. If $M$ is $Spin$ the coefficients of this expansion  
are indices of twisted Dirac operators 
$$\Phi _0(M)=q^{-k/2}\cdot \hat A(M,\bigotimes _{n=2m+1>0}\Lambda _{-q^n}TM \otimes 
\bigotimes 
_{n=2m>0}S_{q^n}TM)$$ 
$$=q^{-k/2}\cdot (\hat A(M) -\hat A(M,TM)\cdot q +\hat A(M,\Lambda ^2{TM}+TM)\cdot 
q^2+\ldots ).$$ Here $\hat A(M,E)$ denotes the index of the Dirac operator twisted with the 
complexification $E_\C $ of a real vector bundle $E$ over $M$. 

For a $Spin$-manifold $M$ the first coefficient of the series $\Phi _0(M)$, the $\hat A$-genus, 
vanishes if $M$ admits a metric of positive scalar curvature \cite{Lic} or if $M$ admits a 
non-trivial smooth $S^1$-action \cite{AtHi}. Our main result asserts that additional coefficients 
vanish if $M$ admits a metric of positive curvature with symmetry rank $\geq 2$. 

\begin{theorem}\label{main theorem} Let $M$ be a closed connected $Spin$-manifold of dimension $> 12r-4$. 
If $M$ admits a metric of positive curvature and symmetry rank $\geq 2r$ then the indices of twisted 
Dirac operators occurring as the first $(r+1)$ coefficients in the expansion $\Phi _0(M)$ vanish. 
\end{theorem}

\noindent
We remark that all simply connected manifolds known to carry a metric of positive curvature have a 
lot of symmetry. Besides the biquotients found by Eschenburg \cite{Es1,Es2} and Bazaikin \cite{Ba}, 
all other examples admit a homogeneous metric of positive curvature. The latter were classified 
by Berger \cite{Be1}, Aloff, Wallach \cite{Wa,AlWa} and B\'erard Bergery \cite{BB} (for recent 
progress on cohomogeneity one manifolds see \cite{GrZi1,GrZi2,SchTu,PoVe2,PoVe1,Ve}). 

In the case that the symmetry rank of $M$ is at least two Theorem \ref{main theorem} states: 
 
\bigskip 
\noindent 
{\it Let $M$ be a closed connected $Spin$-manifold of dimension $> 8$. If $M$ admits a metric of 
positive curvature and symmetry rank $\geq 2$ then $\hat A(M)$ and $\hat A(M,TM)$ vanish.} 
 
\bigskip 
\noindent 
Note that the index $\hat A(M,TM)$ does not vanish for the quaternionic plane. Since the symmetry 
rank of $\H P^2$ (with its standard metric) is three the lower bound on the dimension of $M$ is 
necessary. 

We believe that the vanishing of $\hat A(M,TM)$ also holds under weaker symmetry assumptions. For 
$2$-connected manifolds we show

\begin{theorem}\label{rank one theorem} Let $M$ be a closed $2$-connected manifold of dimension $> 8$. 
If $M$ admits a metric of positive curvature with effective isometric $S^1$-action then $\hat A(M)$ and $\hat A(M,TM)$ vanish. 
\end{theorem}

\noindent
The proofs of Theorem \ref{main theorem} and Theorem \ref{rank one theorem} are rather indirect. For both statements we study the fixed point manifold of isometric cyclic subactions. The Bott-Taubes-Witten rigidity theorem \cite{Wi, Ta,BoTa} for elliptic genera implies that the codimension of the fixed point manifold is bounded from above by a constant which depends on the pole order of the expansion $\Phi _0(M)$ \cite{HiSl,Decy}. Further restrictions arise from the curvature assumption. A component of the fixed point manifold is a totally geodesic submanifold of the positively curved manifold $M$. By an old result of Frankel \cite{Fr1} totally geodesic submanifolds of sufficiently large dimension must intersect. This property imposes additional restrictions on the fixed point manifold. The consequences of the rigidity theorem and Frankel's result indicated above are the main ingredients in the proofs of Theorem \ref{main theorem} and Theorem \ref{rank one theorem} which also rely on the  work of Grove and Searle on isometric $S^1$-actions of codimension two \cite{GrSe1} and recent work of Wilking on the connectivity of the inclusion of totally geodesic submanifolds \cite{Wi}.

We don't know how to prove Theorem \ref{main theorem} and Theorem \ref{rank one theorem} by more 
direct methods such as the Bochner formula for twisted Dirac operators. Already for the proof of the 
vanishing of $\hat A(M,TM)$ we need to use the entire elliptic genus. Note that in view of the above 
discussion for $\H P^2$ a Bochner type argument for the vanishing of $\hat A(M,TM)$ would not apply 
in dimension eight! 

Manifolds of positive curvature (no assumptions on the symmetry) are classified in dimension $<4$ 
\cite{Ham}. In dimension $\geq 4$ the only known obstructions to positive curvature are given by 
restrictions for the fundamental group (cf. \cite{My,Sy,Gr2}, see also \cite{Ro1}), Gromov's Betti 
number theorem \cite{Gr2} and the Lichnerowicz-Hitchin vanishing theorem \cite{Lic,Hit} for the 
$\alpha $-invariant of $Spin$-manifolds. 

Further progress concerning obstructions and classification has been obtained for positively curved 
manifolds with a lot of symmetry, e.g. manifolds with large isometry dimension, large (discrete) 
symmetry rank or small cohomogeneity \cite{HsKl,GrSe1,Wilking,Se,GrSe2,Ro2,SeYa,Ya,PuSe,FaRo,Deinprep}. 

All these results require that the dimension of the manifold is bounded from above by a constant 
depending on the symmetry. In contrast Theorem \ref{main theorem} and Theorem \ref{rank one theorem} 
only require a {\it lower} bound on the dimension of the manifold. 

Theorem \ref{main theorem} allows to distinguish positive curvature from weaker curvature properties 
under assumptions on the symmetry rank. For example, consider the product of $\H P^2$ and a 
Ricci-flat $K_3$-surface. The Riemannian manifold $M=\H P^2\times K_3$ (equipped with the product 
metric) has symmetry rank three and positive scalar curvature as well as non-negative Ricci 
curvature. The index $\hat A(M,TM )$ does not vanish. Hence, if one restricts to metrics with 
symmetry rank $\geq 2$ it follows from Theorem \ref{main theorem} that $M$ admits a metric of 
positive scalar curvature but no metric of positive curvature. 

This kind of reasoning can be pushed further to yield examples of simply connected manifolds of 
positive Ricci curvature which do not admit a metric of positive curvature if one restrict to metrics 
with a prescribed lower bound on the symmetry rank (see Section \ref{versus} for precise statements). 
Using different arguments (based on \cite{Decy,Wilking}) it is possible to distinguish positive Ricci 
from positive curvature under rather mild assumptions on the symmetry. We shall call an $S^1$-action 
on a Riemannian manifold {\it \almostisometric } of order $\oorder $ if the cyclic subgroup of order 
$\oorder $ acts effectively and isometrically. 
\pagebreak
\begin{theorem}\label{small Betti} For every $d\in \N $ and every $\oorder \geq 2 $ there exists a 
simply connected closed manifold $M$ of dimension greater than $d$ such that: 
\begin{enumerate} 
\item $M$ admits a metric of positive Ricci curvature with \almostisometric \ $S^1$-action of order $\oorder$. 
\item $M$ does not admit a metric of positive curvature with \almostisometric \ $S^1$-action of order 
$\oorder$. \end{enumerate} 
\end{theorem}

\noindent
We note that the examples (given in Section \ref{versus}) admit metrics of positive Ricci curvature 
and symmetry rank $\geq 3$. In particular, one obtains simply connected manifolds of positive Ricci 
curvature and positive symmetry rank which do not admit a metric of positive curvature with positive 
symmetry rank. 

The paper is structured in the following way. In the next section we review basic properties of 
positive curvature used in the proofs of Theorem \ref{main theorem}, Theorem \ref{rank one theorem} 
and Theorem \ref{small Betti}. In Section \ref{rigidity} we recall the rigidity theorem for elliptic 
genera and discuss applications to cyclic actions. The proofs of Theorem \ref{main theorem} and 
Theorem \ref{rank one theorem} are given in Section \ref{proof of main theorem} and Section 
\ref{proof of main lemma}. In Section \ref{proof of main theorem} we also discuss related results for 
positive $k$th Ricci curvature, finite isometric actions and integral cohomology $\H P^k$'s. 
 In the final section we 
show Theorem \ref{small Betti}. 

\section{Geodesic submanifolds}\label{geodesic}
In this section we review basic properties of positively curved manifolds used in the proofs of 
Theorem \ref{main theorem}, Theorem \ref{rank one theorem} and Theorem \ref{small Betti}. A main 
ingredient is an old result of Frankel on the intersection property for totally geodesic 
submanifolds. 

\begin{theorem}[ \cite{Fr1}]
\label{intersection theorem} Let $N_1$ and $N_2$ be totally ge\-odesic 
submanifolds of a positively curved connected manifold $M$. If $\dim N_1+ 
\dim N_2\geq \dim M$ then $N_1$ and $N_2$ have non-empty intersection.\proofend
\end{theorem}

\noindent The proof uses a Synge type argument for the parallel transport along a geodesic 
from $N_1$ to $N_2$ which minimizes the distance. Whereas it is difficult to find totally geodesic 
submanifolds for generic metrics they do occur naturally as fixed point components in the presence of 
symmetry. Theorem \ref{intersection theorem} clearly imposes restrictions on the fixed point manifold 
of isometric actions. The following consequence is immediate. 

\begin{corollary}\label{frankel corollary} Let $\sigma $ be an isometry of a positively curved connected manifold 
$M$ and let $F$ be a connected component of the fixed point manifold $M^\sigma$ of minimal 
codimension. Then the dimension of every other component is less than the codimension of 
$F$.\proofend 
\end{corollary}

\noindent
In \cite{Fr2} Frankel applied Theorem \ref{intersection theorem} to show that the inclusion $N\hra M$ 
of a totally geodesic submanifold of codimension $k$ is $1$-connected provided $k$ is not larger than 
half of the dimension of the positively curved manifold $M$. Recently, Wilking generalized this 
result significantly using Morse theory. 
\begin{theorem}[\cite{Wilking}]\label{connectivity theorem} Let $M$ be
 an $n$-dimensional connected
 Riemannian manifold and $N$ a connected totally geodesic submanifold 
of codimension $k$. If $M$ is positively curved then the inclusion $N\hra M$ is 
$(n-2k+1)$-connected.\proofend 
\end{theorem}

\noindent
The theorem imposes severe restrictions on the topology. For example, if $M$ is simply connected and 
the codimension of $N$ is two then either $M$ is homeomorphic to an odd dimensional sphere or all odd 
Betti numbers of $M$ and $N$ vanish \cite{Wilking}. The case where $N$ is fixed under an isometric 
$S^1$-action was studied before by Grove and Searle. 

\begin{theorem}[\cite{GrSe1}]\label{small codim theorem} Let $M$ be a simply connected positively curved 
manifold with isometric $S^1$-action. If $\codim M^{S^1}=2$ then $M$ is diffeomorphic to a sphere or 
a complex projective space.\proofend 
\end{theorem}

\section{Rigidity and cyclic actions}\label{rigidity}
In this section we recall the rigidity theorem for elliptic genera and discuss applications to cyclic 
actions. For more information on elliptic genera we refer to \cite{La,HiBeJu}. 

A genus is a ring homomorphism from the oriented bordism ring $\Omega _*^{SO}$ to a $\Q $-algebra $R$ 
\cite{Hi1}. The genus is called elliptic (of level $2$) if its logarithm $g(u)$ is given by a formal 
elliptic integral 
$$g(u)=\int _0^u\frac {dt}{\sqrt {1-2\cdot \delta \cdot t^2+\epsilon \cdot t^4}},$$
where $\delta , \epsilon\in R$ \cite{Oc}. Classical examples of elliptic genera are the signature 
($\delta 
=\epsilon =1$) and the $\hat A$-genus ($\delta =-\frac 1 8$, $\epsilon 
=0$). 

The ring of modular forms $M_*(\Gamma 
_0(2))$ is a polynomial ring with generators $\delta $ and $\epsilon$ of
weight $2$ and $4$, respectively \cite{HiBeJu}. The corresponding elliptic genus 
$$\varphi :\Omega 
_*^{SO}\to M_*(\Gamma 
_0(2))$$
is universal since $\delta $ and $\epsilon $ are algebraically independent. 

As in \cite{Wi,HiSl} we shall consider for a $4k$-dimensional oriented manifold $M$ the normalized 
elliptic genus $\Phi (M)=\varphi (M)/\epsilon ^{k/2}$ which is a modular function of weight $0$ (with 
$\Zp 2$-character). In one of the cusps (the signature cusp) $\Phi (M)$ has an expansion which is 
equal to the following series of twisted signatures 
$$sign(q,{\cal L}M):=sign (M,\bigotimes _{n=1}^\infty 
S_{q^n}TM 
\otimes 
\bigotimes 
_{n=1}^\infty \Lambda_{q^n}TM  )\in \Z \lbrack\lbrack q\rbrack \rbrack .$$
Here $sign(M,E)$ denotes the index of the signature operator twisted with the complexified \vb \ 
$E_\C $. The series $sign(q,{\cal L}M)$ describes the ``signature'' of the 
free loop space ${\cal L}M$ localized at the manifold $M$ of constant loops \cite{Wi}. 

In the other cusp (the $\hat A$-cusp) $\Phi (M)$ expands as a series of characteristic 
numbers\footnote{For a suitable change of cusps.}
\pagebreak
$$\Phi _0(M):=q^{-k/2}\cdot \hat A(M,\bigotimes _{n=2m+1>0}\Lambda _{-q^n}TM \otimes 
\bigotimes 
_{n=2m>0}S_{q^n}TM)$$ 
$$=q^{-k/2}\cdot (\hat A(M) -\hat A(M,TM)\cdot q +\hat A(M,\Lambda ^2{TM}+TM)\cdot 
q^2+\ldots ).$$ Here $\hat A(M,E):=\langle \hat {\cal A}(M)\cdot ch(E_\C), \lbrack M\rbrack \rangle 
$, where $ \hat {\cal A}(M)$ denotes the 
multiplicative sequence for the $\hat A$-genus, $E_\C $ is the complexification of the vector bundle $E$, $\lbrack M\rbrack $ denotes the fundamental cycle and 
$\langle \quad, \quad \rangle $ is the Kronecker pairing. If $M$ is $Spin$ $\hat A(M,E)$ is equal to 
the index of the Dirac operator twisted with $E_\C $ by the Atiyah-Singer index theorem 
\cite{AtSiIII}. In this case $\Phi _0(M)$ has an interpretation as a series of indices of twisted 
Dirac operators (twisted Dirac-indices for short).

The main feature of the elliptic genus is its rigidity under actions of compact connected Lie groups. 
The rigidity was explained by Witten in \cite{Wi} using heuristic arguments from quantum field theory 
and proved rigorously by Taubes and Bott-Taubes \cite{Ta,BoTa} (cf. also \cite{Hi3, Li}). 

Assume the $Spin$-manifold $M$ carries an action by a compact Lie group $G$ preserving the 
$Spin$-structure (note that any smooth $G$-action lifts to the $Spin$-structure after passing to a two-fold covering action, if necessary). Then the indices occurring in the expansions of $\Phi (M)$ refine to virtual 
$G$-representations which we identify with their characters. If $G$ is connected the elliptic 
genus is rigid, i.e. the characters do not depend on $G$. 

\begin{theorem}[\cite{Ta,BoTa}]\label{rigidity theorem} Let $M$ be a $G$-equivariant $Spin$-manifold.
 If $G$ is 
connected 
 then each twisted signature (resp. each twisted Dirac-index)
 occurring as coefficient in the expansion of $\Phi (M)$ in the signature cusp (resp. in the $\hat A$-cusp) is constant as a character of $G$.\proofend 
\end{theorem}

\noindent
The rigidity of $\Phi (M)$ also holds for certain non-$Spin$ manifolds such as $Spin^c$-manifolds 
with first Chern class a torsion class \cite{Despinc} or orientable manifolds with finite second 
homotopy group \cite{herrera}.

In the remaining part of this section we discuss consequences of Theorem \ref{rigidity theorem} for 
cyclic actions which are used in the proofs of Theorem \ref{main theorem}, Theorem \ref{rank one 
theorem} and Theorem \ref{small Betti}. 

Assume $S^1$ acts on the $Spin$-manifold $M$ (not necessarily preserving the $Spin$-structure). Let 
$\sigma \in S^1$ be the element of order $2$. In \cite{HiSl} Hirzebruch and Slodowy showed that the 
expansion of the elliptic genus in the signature cusp can be expressed in terms of the transversal 
self-intersection $M^\sigma \circ M^\sigma $ of the fixed point manifold $M 
^\sigma $
\begin{equation}\label{selfinter}sign 
(q,{\cal L}M)=sign(q,{\cal L}(M^\sigma \circ M^\sigma )).\end{equation} The formula is equivalent to 
$\Phi _0(M)=\Phi _0(M^\sigma \circ M^\sigma )$ which implies the following generalization of the 
Atiyah-Hirzebruch vanishing theorem \cite{AtHi} for the $\hat A$-genus. 
\begin{theorem}[\cite{HiSl}]\label{involution} Let $M$ be a $Spin$-manifold with 
$S^1$-action and let 
 $\sigma $ be the element of order two in $S^1$. If $\codim M^\sigma > 4r$
 then the first $(r+1)$ coefficients of $\Phi _0(M)$ vanish.\proofend
\end{theorem}

\noindent Here
$\codim M^\sigma $ denotes the minimal codimension of the connected components of $M^\sigma $ in $M$. 

Recall that the $S^1$-action is called even if it lifts to the $Spin$-structure (otherwise the action 
is called odd). In the even case the codimension of all fixed point components of $M^\sigma $ is 
divisible by $4$ whereas in the odd case the codimensions are always $\equiv 2 \bmod 4$ (cf. 
\cite{AtHi}, Lemma 2.4). Note that for an odd action the series $q^{\dim M/8}\cdot \Phi _0(M)$ is an 
element in $\C \lbrack \lbrack q\rbrack 
\rbrack $ whereas $q^{\dim M/8}\cdot 
\Phi _0(M^\sigma \circ M^\sigma )\in q^{1/2}\cdot \C 
\lbrack \lbrack q\rbrack \rbrack $. Thus formula (\ref{selfinter}) implies

\begin{corollary}\label{oddaction} Let $M$ be a $Spin$-manifold with $S^1$-action. 
If the action is odd then $\ell (M)$ vanishes identically.\proofend 
\end{corollary}

\noindent
In the remaining part of this section we recall a generalization of Theorem \ref{involution} to 
cyclic actions of arbitrary order \cite{Decy}. Let $M$ be a $Spin$-manifold with $S^1$-action and let 
$\sigma \in S^1$ be of order $\oorder \geq 2$. At a connected component $Y$ of the fixed point 
manifold $M^{S^1}$ the tangent bundle $TM$ splits equivariantly as the direct sum of $TY$ and the 
normal bundle $\nu $. The latter splits (non-canonically) as a direct sum $\nu 
=\bigoplus_{k\neq 0} \nu _k$ corresponding to the irreducible real $2$-dimensional 
$S^1$-representations $e^{i\cdot 
\theta }\mapsto \left (\begin{smallmatrix} 
\cos k\theta & \sin k\theta \\ 
-\sin k\theta &\cos k\theta
\end{smallmatrix}\right )$, $k\neq 0$. We fix such a decomposition of $\nu $. 
For each $k\neq 0$ choose $\alpha 
_k\in \{\pm 1\}$ such that $\alpha _k k\equiv \tilde k \bmod \oorder $, 
$\tilde k\in 
\{0,\ldots , \lbrack 
\frac o 2\rbrack \}$.
 On each vector bundle 
$\nu _k$ introduce a complex structure such that $\lambda \in S^1$ acts on $\nu _k$ by scalar 
multiplication with $\lambda ^{\alpha _k k}$. Finally define 
$$m_\oorder (Y):=(\sum _k \dim  _\C  \nu_k 
\cdot \tilde k )/\oorder \quad \text{ and }\quad m_\oorder :=\min _{Y} m_\oorder (Y),$$
where $Y$ runs over the connected components of $M^{S^1}$. 
\begin{theorem}[\cite{Decy}]\label{theorem cyclic} Let $M$ be a $Spin$-manifold with 
$S^1$-action. If $m_\oorder > r$ then the first $(r+1)$ coefficients of $\Phi _0(M)$ vanish.\proofend 
\end{theorem}

\noindent
To prove this result one analysis the expansion of the equivariant elliptic genus in the $\hat 
A$-cusp using the \LFF \ \cite{AtSiIII} and Theorem \ref{rigidity theorem} (see \cite{Decy} for details). Since 
$\codim M^\sigma 
\leq 2\oorder \cdot m_\oorder $ Theorem \ref{theorem cyclic} implies the following 
corollary which generalizes Theorem \ref{involution} to cyclic actions of arbitrary finite order. 
 
\begin{corollary}\label{corollary special case}
 Let $M$ be a $Spin$-manifold with 
$S^1$-action and let $\sigma \in S^1$ be of order $\oorder 
\geq 2$. If $\codim M^\sigma 
> 2\oorder \cdot r$ then the first $(r+1)$ coefficients of $\Phi _0(M)$ vanish.\proofend 
\end{corollary}

\section{Positive curvature and elliptic genera}\label{proof of main theorem}
In this section we prove the vanishing results for positively curved manifolds stated in the 
introduction. Recall from Section \ref{geodesic} that positive curvature restricts the 
dimension of fixed point components of isometric actions. This property is an essential ingredient in 
the proofs of Theorem \ref{main theorem} and Theorem \ref{rank one theorem}. 

\begin{definition} {\rm Let $M$ be a closed connected
 manifold with smooth action by a torus $T$. We say $M$
 has {\it \rfpd }\ for the prime $p$ if for all cyclic subgroups $H\subset H^\prime \subset T$ of order a power of $p$ and every connected component $X\subset M^H$ the dimension of two different connected components $F_1$ and $F_2$ of $X^{H^\prime}$ is restricted by $\dim F_1 + \dim F_2 <\dim X$.}
\end{definition}

\noindent
By Theorem \ref{intersection theorem} a positively curved manifold with isometric $T$-action has 
\rfpd \ for any prime $p$. Other examples are given by $T$-manifolds with the same integral cohomology ring as a 
projective space (see \cite{Bredon}, Chapter VII). 

A main step in the study of the elliptic genus for these manifolds is the following lemma which is a consequence of Theorem \ref{involution} and Theorem 
\ref{theorem cyclic} (see the next section for the proof). 

\begin{lemma}[Main lemma]\label{mlemma}
Let $T$ be a torus of rank $2r>0$ and let $M$ be a $Spin$-manifold with effective $T$-action. If $M$ 
has \rfpd \ for the prime $2$ and $\dim M>12r-4$ then at least one of the following possibilities holds: 
\begin{enumerate}
\item The first $(r+1)$ coefficients in the expansion $\Phi _0(M)$ vanish.
\item For some subgroup $S^1\subset T$ with involution $\sigma \in S^1$ the fixed point manifolds $M^\sigma $ and $M^{S^1}$ are orientable and connected of codimension $4$ and $6$, respectively.
\end{enumerate}
\end{lemma}

\noindent
With this information at hand we now prove Theorem \ref{main theorem} stated in the introduction. 

\begin{theorem}\label{main theorem I} Let $M$ be a closed connected $Spin$-manifold of dimension $> 12r-4$. 
If $M$ admits a metric of positive curvature and symmetry rank $\geq 2r$ then the indices of twisted 
Dirac operators occurring as the first $(r+1)$ coefficients in the expansion $\Phi _0(M)$ vanish. 
\end{theorem}

\noindent
{\bf Proof:} For $r=0$ the theorem follows from \cite{Lic}. So assume $r>0$. Let $T$ denote a torus 
of rank $2r$ acting isometrically and effectively on $M$. Assume the first $(r+1)$ coefficients in 
the expansion $\Phi _0(M)$ do not vanish. In particular, $M$ is of even dimension and simply 
connected \cite{Sy}. 

The $T$-action has the properties given in the second part of Lemma \ref{mlemma}, i.e. $X:=M^\sigma $ 
and $X^{S^1}$ are orientable connected submanifolds of codimension $4$ and $6$, respectively. Note 
that $X$, being a totally geodesic submanifold of $M,$ inherits positive curvature from $M$. Since 
the submanifold $X^{S^1}$ of $X$ has codimension two it follows from Theorem \ref{small codim 
theorem} that $X$ is diffeomorphic to a sphere or a complex projective space. It is well-known that a 
non-trivial $S^1$-action on a complex projective space has more than one connected fixed point 
component (cf. \cite{Bredon}, Chapter VII, Theorem 5.1). Thus $X$ is diffeomorphic to a sphere of codimension $4$. Since $\dim M>8$ 
the Euler class $e(\nu _X)\in H^4(X;\Z )$ of the normal bundle $\nu _X$ of $X\subset M$ is trivial 
(we fix compatible orientations for $X$ and $\nu _X$). 

This implies that the expansion of the elliptic genus in the signature cusp 
vanishes as we will explain next (see \cite{Decy} for details). Recall that the coefficients of this expansion are twisted signatures $sign (q,{\cal L}M)=\sum _{l\geq 0}sign (M,E_l)\cdot q^l$, where $E_l$ is a virtual complex \vb \ associated to the tangent bundle $TM$. The $S^1$-action on $TM$ induces an action on each $E_l$. Let $sign _{S^1}(M,E_l)$ denote the $S^1$-equivariant twisted signature and $sign _{S^1}(q,{\cal L}M)$ the equivariant expansion of the elliptic genus. Recall from Theorem \ref{rigidity theorem} that $sign (q,{\cal L}M)=sign _{S^1}(q,{\cal L}M)(\lambda )$ for any $\lambda \in S^1$. By the \LFF \ \cite{AtSiIII} $sign _{S^1}(M,E_l)(\sigma )$ is equal to 
$$\left \langle \prod _i\left (x_i\cdot \frac {1+e^{-x_i}}{1-e^{-x_i}}\right )\cdot \prod _j
\left (y_j\cdot\frac {1+e^{-y_j}}{1-e^{-y_j}}\right )^{-1}
\cdot ch (E_l)(\sigma )\cdot e(\nu _X), \lbrack X\rbrack \right \rangle .$$
Here $\pm x_i$ (resp. $\pm y_j$) denote the formal roots of $X$ (resp. the normal bundle $\nu _X$), 
$ch(E_l)$ denotes the equivariant Chern character of $E_l$, $\lbrack X\rbrack $ the fundamental cycle 
and $\langle 
\quad ,\quad \rangle $ the 
Kronecker pairing. 

Since $e(\nu _X)$ is trivial $sign _{S^1}(M,E_l)(\sigma )=0$ for every $l\geq 0$. By Theorem 
\ref{rigidity theorem} $sign (q,{\cal L}M)$ vanishes identically. 
Thus $\Phi_0 (M)=0$ contradicting the initial assumption of the proof.\proofend 

\begin{remarks}\begin{enumerate}
\item Theorem \ref{main theorem I} also follows from Lemma \ref{mlemma} and Theorem \ref{connectivity theorem}. 
\item The statement for $r=2$ in Theorem \ref{main theorem I} can be strengthened to: If $M$ admits a metric of positive curvature
  and symmetry rank $\geq 3$ then the first three coefficients in the expansion $\Phi _0(M)$ vanish,
   i.e. $\hat A(M)=\hat A(M, TM)=\hat A(M,\Lambda^2(TM))=0$.
\end{enumerate}
\end{remarks}

\noindent
Under stronger assumptions on the bounds for the symmetry rank and the dimension the conclusion of 
Theorem \ref{main theorem I} holds if one only assumes that $M$ has positive $k$th Ricci curvature 
and an elementary $p$-abelian subgroup of the torus acts by isometries.

Recall that a manifold $M$ 
has positive $k$th Ricci curvature (or $k$-positive Ricci curvature) if for any $(k+1)$ mutually 
orthogonal unit tangent vectors $e, e_1,\ldots ,e_k$ (at any point of $M$) the sum of curvatures 
$\sum _{i=1}^k sec(e\wedge e_i)$ is positive \cite{Ken}. Thus, $1$-positive Ricci curvature is 
equivalent to positive curvature and $(\dim M-1)$-positive Ricci curvature is equivalent to positive 
Ricci curvature.

Assume that $M$ has positive $k$th Ricci curvature and assume that a torus $T$ of 
rank $R$ acts smoothly on $M$ such that the induced action of the $p$-torus $T_p\cong (\Zp p)^R$, $p$ 
a prime, is isometric and effective. To keep the exposition simple we shall assume the generous 
bounds $R\geq p\cdot r +\frac {k+1} 2$ and $\dim M\geq 6p\cdot r +(k-1)$. 

\begin{proposition}\label{kth ricci} For a connected $Spin$-manifold $M$ as above the indices of twisted 
Dirac operators occurring as the first $(r+1)$ coefficients in the expansion $\Phi _0(M)$ vanish. 
\end{proposition}

\noindent
{\bf Sketch of proof:} First note that the intersection property for totally geodesic submanifolds in 
positive curvature (Theorem \ref{intersection theorem}) extends to positive $k$th Ricci curvature 
\cite{Ken}: {\it Two totally geodesic submani\-folds $N_1$ and $N_2$ of a manifold $M$ of positive 
$k$th Ricci curvature intersect if $\dim N_1+\dim N_2\geq \dim M+(k-1)$.} 

In particular, if $F_1$ and $F_2$ are two different connected fixed point components of an isometry 
$\sigma $ then $\dim F_1 + \dim F_2 <\dim M+(k-1)$. 

Assume the first $(r+1)$ coefficients in the expansion $\Phi _0(M)$ do not vanish. Consider the 
action of $\sigma \in T_p$ on $M$. By Corollary \ref{corollary special case} the codimension of $M^\sigma 
$ is $\leq 2 p\cdot r$. Hence, a connected component $F$ of $M^\sigma $ has either ``small 
codimension'', i.e. $\codim F\leq 2 p\cdot r$, or ``small dimension'', i.e. $\dim F < 2 p\cdot r+(k-1)$.

Consider a $T$-fixed point $pt\in M$ (which exists since $\Phi (M)\neq 0$) and let $F_\sigma 
\subset M^\sigma$ denote the component which contains $pt$. It is an elementary exercise to show that 
the $p$-torus $T_p$ has a basis $\sigma 
_1,\ldots ,\sigma _R$ such that $\dim F_{\sigma _i}\geq 2R-2\geq 2p\cdot r +(k-1)$. This implies that the 
codimension of $F_{\sigma _i}$ is small. Since $\dim M\geq 6p\cdot r +(k-1)$ the codimension of 
$F_{\sigma +\sigma ^\prime }$ is small provided this holds for $F_{\sigma}$ and $F_{\sigma^\prime }$. 
Hence, the codimension of $F_{\sigma}$ is small for every $\sigma \in T_p$, i.e. $\codim F_\sigma 
\leq 2 p\cdot r$ for every $\sigma \in T_p$. However, it follows from elementary linear algebra that for some $\sigma \in T_p$ 
the codimension of $F_\sigma $ is at least $2R\geq 2p\cdot r +(k+1)$. This gives the desired 
contradiction. Hence, the first $(r+1)$ coefficients in the expansion $\Phi _0(M)$ vanish.\proofend 

\bigskip
\noindent
The next result implies Theorem \ref{rank one theorem} stated in the introduction. 
\begin{theorem}\label{rank one theorem I} Let $M$ be a closed $2$-connected manifold of dimension $> 8$ with effective smooth $S^1$-action. 
Assume $M$ admits a metric of positive curvature such that the 
subgroup $\Zp 4\subset S^1$ acts by isometries. Then $\hat A(M)$ and $\hat A(M,TM)$ vanish. 
\end{theorem}

\noindent
{\bf Proof:} We fix an orientation and a $Spin$-structure for $M$. Note that $\hat A(M)$ vanishes by \cite{Lic,AtHi}. Let 
$\rho \in S^1$ be an element of order $4$ and let $\sigma :=\rho ^2$. Assume $\hat A(M,TM)$ does not 
vanish. Then $M$ has dimension $4k\geq 12$ and the fixed point manifold $M^\sigma $ is the union of a 
connected component $X$ of codimension $4$ and (a possible empty set of) isolated fixed points by Corollary \ref{frankel 
corollary}, Theorem \ref{involution} and Corollary \ref{oddaction}. Next apply Theorem \ref{theorem 
cyclic} (for $\oorder =4$) to conclude that either $X=X^\rho$ or the codimension of $X^\rho $ in $X$ 
is two. 

We claim that $X=X^\rho$. If the codimension of $X^\rho$ is two then the inclusion $X^\rho \hra X$ is 
$(4k-7)$-connected by Theorem \ref{connectivity theorem}. This implies that $b_4(X)$ is equal to 
$b_2(X)$ (use Poincar\'e duality for $X$ and $X^\rho $). Note that $X$ is $2$-connected by Theorem 
\ref{connectivity theorem}. Hence, $b_4(X)$ vanishes. In particular, the Euler class of the normal 
bundle of $X$ in $M$ vanishes rationally. Now argue as in the proof of Theorem \ref{main theorem I} 
to conclude that $sign(q,{\cal L}M)=sign_{S^1}(q,{\cal L}M)(\sigma )=0$ contradicting the assumption on 
$\hat A(M,TM)$. Hence, $X=X^\rho$. 

Since $X$ is fixed by $\rho $ the action of $\rho $ on the normal bundle $\nu _X$ induces a complex 
structure such that $\rho $ acts by multiplication with $i$. We fix the orientation of $X$ which is 
compatible with the orientations of $\nu _X$ and $M$. Also the action of $\rho $ induces a complex 
structure on the normal bundle of any of the isolated $\sigma $-fixed points. We shall now compute 
the local contributions in the 
\LFF \ for 
$sign_{S^1}(q,{\cal L}M)(\rho )$: 
$$sign_{S^1}(q,{\cal L}M)(\rho )=\mu _X +\sum \mu _{p_j}$$
Here $\mu _X$ (resp. $\mu _{p_j}$) denotes the local contribution at $X$ (resp. at an isolated fixed 
point $p_j$). The term $\mu _X$ is given by \cite{AtSiIII} 
$$\mu _X=\left \langle {\cal T}_X\cdot {\cal N}_X, \lbrack X\rbrack \right \rangle $$
where 
$${\cal T}_X:= \prod _i\left (x_i\cdot \frac {1+e^{-x_i}}{1-e^{-x_i}}\cdot 
\prod _{n=1}^\infty \frac {(1+ q^n\cdot e^{x_i})\cdot (1+ q^n\cdot e^{-x_i})}{(1- q^n\cdot e^{x_i})
\cdot (1- q^n\cdot e^{-x_i})}\right )$$
and 
$${\cal N}_X:= \prod _{j=1,2} \left (\frac {1+e^{-y_j}\cdot (-i)}{1-e^{-y_j}\cdot (-i)}\cdot 
\prod _{n=1}^\infty \frac {(1+ q^n\cdot e^{y_j}\cdot i)\cdot (1+ q^n\cdot e^{-y_j}\cdot (-i))}
{(1- q^n\cdot e^{y_j}\cdot i)\cdot (1- q^n\cdot e^{-y_j}\cdot (-i))}\right )$$ Here $\pm x_i$ (resp. 
$y_1, y_2$) denote the formal roots of $X$ (resp. $\nu _X$). Since $X$ is $2$-connected the first 
Chern class $y_1+y_2$ of $\nu _X$ vanishes. This implies 
$${\cal N}_X=\frac {1+e^{-y_1}\cdot (-i)}{1-e^{-y_1}\cdot (-i)}\cdot \frac {1+e^{y_1}
\cdot (-i)}{1-e^{y_1}\cdot (-i)}=-1.$$
Hence, the expression for $\mu _X$ simplifies to 
$$\mu _X=-\left \langle {\cal T}_X, \lbrack X\rbrack \right \rangle =-sign(q,{\cal L}X).$$
The term $\mu _{p_j}$ is given by \cite{AtSiIII} 
$$\mu _{p_j}=\pm \left (\frac {1-i}{1+i}\prod _{n=1}^\infty \frac {(1+i\cdot q^n)
\cdot (1-i\cdot q^n)}{(1-i\cdot q^n)\cdot (1+i\cdot q^n)}\right )^{2k}=
\pm  \left (\frac {1-i}{1+i}\right )^{2k}=\pm 1$$
By Theorem \ref{rigidity theorem} $sign(q,{\cal L}M)=sign_{S^1}(q,{\cal L}M)(\rho )=-sign(q,{\cal L}X)+c$, where $c$ is the integer obtained by summing up $\mu _{p_j}$. Equivalently, 
\begin{equation}\label{rank one equation} \Phi(M)=-\Phi(X)+c\end{equation}
Note that $\Phi_0(M)\in q^{-k/2}\C 
\lbrack \lbrack q\rbrack \rbrack $, whereas $\Phi_0(X)\in q^{1 /2}\cdot q^{-k/2}\C \lbrack \lbrack q\rbrack 
\rbrack $. By comparing the expansions in the $\hat A$-cusp of both sides of (\ref{rank one 
equation}) it follows that $\Phi_0(M)\in \Z $ (in fact, depending on the parity of $k$, $\Phi_0(M)$ 
is equal to $0$ or $c$). Since $\dim M>8$ this implies that $\hat A(M,TM)=0$.\proofend

\bigskip
\noindent
In the remaining part of this section we use Lemma \ref{mlemma} to study the Pontrjagin numbers of 
$4k$-dimensional $Spin$-manifolds with symmetry which have the same integral cohomology ring as a 
projective space. Such manifolds are either integral cohomology $\H P^k$'s or integral cohomology 
Cayley planes (for recent progress in the study of integral cohomology $\C P^m$'s with symmetry we 
refer to \cite{Decpn} and references therein). Note that an integral cohomology Cayley plane with 
smooth non-trivial $S^1$-action has the same Pontrjagin numbers as the Cayley plane since these are 
completely determined by the signature and the $\hat A$-genus (the latter vanishes by \cite{AtHi}). 
The same argument applies to an integral cohomology $\H P^2$ with smooth non-trivial $S^1$-action. We 
shall now apply Lemma \ref{mlemma} to integral cohomology $\H P^k$'s for $k>2$. 

\begin{proposition}\label{hpn proposition} Let $M$ be an integral cohomology $\H P^k $, i.e. $H^*(M;\Z 
)\cong \Z \lbrack x\rbrack /(x^{k+1})$, where $x$ has degree $4$. If a torus $T$ of rank $2r>0$ acts 
effectively and smoothly on $M$ and $4k>12r-4$ then the first $(r+1)$ coefficients in the expansion 
$\Phi _0(M)$ vanish. 
\end{proposition}

\noindent
{\bf Proof:} Note that $M$ is a $Spin$-manifold with \rfpd \ for any prime. By Lemma \ref{mlemma} we may 
assume that for some $S^1$-subgroup of $T$ the fixed point manifold $M^{S^1}$ is connected of 
codimension $6$. It is well known that $M^{S^1}$ is an integral cohomology projective space of the 
form $H^*(M^{S^1};\Z )\cong \Z 
\lbrack y\rbrack /(y^{l+1})$ where $y$ has even degree $\leq 4$ (cf. \cite{Bredon}, Chapter VII, Theorem 5.1). Also $\chi (M)=\chi 
(M^{S^1})$ by the \LFF \ for the Euler characteristic. Since the codimension of $M^{S^1}$ is equal to 
$6$ it follows that $\deg (y)=2$ and $k=l=3$. Hence, $M$ is an integral cohomology $\H P^3$. For a 
$12$-dimensional manifold the elliptic genus is a linear combination of the signature and the $\hat 
A$-genus. The signature of $M$ vanishes for trivial reasons. The $\hat A$-genus of $M$ vanishes since 
$M$ is a $Spin$-manifold with non-trivial $S^1$-action \cite{AtHi}. Hence, $\Phi (M)=0$. This 
completes the proof.\proofend 

\noindent
Finally we point out the following consequence of Lemma \ref{mlemma}. 

\begin{proposition}
Let $T$ be a torus of rank $2r>0$ and let $M$ be a $Spin$-manifold with $T$-action of dimension 
$>12r-4$. Assume the signature of $M$ does not vanish. If $M$ has \rfpd \ for the prime $2$ then the first $(r+1)$ 
coefficients in the expansion $\Phi _0(M)$ vanish. 
\end{proposition}

\noindent
{\bf Proof:} By the rigidity of the signature $sign (M)=sign (M^{S^1})$ for any $S^1$-action. Since the 
signature of $M$ does not vanish $\dim M\equiv \dim M^{S^1}\equiv 0 \bmod 4$ if $M^{S^1}$ is connected. In particular, 
$\codim M^{S^1}\neq 6$ and the statement follows from Lemma \ref{mlemma}.\proofend 

\section{Proof of the main Lemma}\label{proof of main lemma}
In this section we prove Lemma \ref{mlemma}. Let $T$ be a torus of rank $2r>0$ and let $M$ be a 
$Spin$-manifold with effective $T$-action. Assume $M$ has \rfpd \ for the prime $2$ and $\dim M>12r-4$. Assume also 
that the first $(r+1)$ coefficients in the expansion $\Phi _0(M)$ do not vanish. This implies that 
$\dim M=4k\geq 12r$.

Our goal is to show that for some subgroup $S^1\subset T$ with involution $\sigma \in S^1$ the fixed point manifolds $M^\sigma $ and $M^{S^1}$ are orientable and connected of codimension $4$ and $6$, respectively. To this end we first examine the action of $T$ at a fixed point $pt$ 
(which exists since $\Phi (M)\neq 0$). 

The tangent space $T_{pt} M$ of $M$ at $pt$ decomposes (non-canonically) into $2k$ complex 
one-dimensional $T$-representations $T_{pt}M=R_1\oplus \ldots \oplus R_{2k}$. With res\-pect to such 
a decomposition the action of $T$ on $T_{pt}M$ is given by a homomorphism $T\hra U(1)^{2k}$. We 
denote by $h$ the induced homomorphism of integral lattices $I_T\to  I_{U(1)^{2k}}=\Z ^{2k}$. Since the $T$-action is effective the mod $p$ reduction
 $I_T\otimes (\Zp p) \to (\Zp p)^{2k}$ of $h$ is injective for every 
$p\geq 2$.

For a finite covering homomorphism $\hat T:=S^1\times 
\ldots \times S^1\to T$ we denote by $\hat h$ the homomorphism 
of integral lattices $I_{\hat T}=\Z ^{2r}\to 
\Z ^{2k}$ for the covering action on $T_{pt} M$. Note that an element $\hat h(a)$, $a\in I_{\hat T}$, determines a homomorphism $S^1\overset \iota \to \hat T$, $\exp (\theta )\mapsto \exp (\theta \cdot a)$, and an $S^1$-subgroup $\iota (S^1)\subset \hat T$. Let $c_1,\ldots ,c_{2r}$ be the image of the standard basis of $I_{\hat T}=\Z 
^{2r}$ under $\hat h$.

\noindent
\begin{claim}\label{first claim} Let $p$ be a prime. There exists a finite covering homomorphism 
$\hat T:=S^1\times 
\ldots \times S^1\to T$ of degree coprime to $p$
 such that
$$A:=\left (\begin{smallmatrix}c_1\\ 
\vdots \\ c_{2r}\end{smallmatrix}\right ) =\left (\begin{smallmatrix}? & 0 & 0 & 0 & 
\ldots & 0 & ? & ? & ? & \ldots & ?\\ 0 & ? & 0 & 0 & \ldots & 0 & ? & 
? & ? & \ldots & ?\\ 0 & 0 & ?& 0 & 
\ldots & 0 & ? & ? & ? & \ldots & ?\\ 
\vdots &\vdots &\vdots & \vdots & \ddots & \vdots  & \vdots & \vdots & \vdots & \vdots & \vdots \\
0&0 &0 &0 &\ldots & ?& ? &? &? &\ldots &? \\ 
\end{smallmatrix}\right ) 
\equiv 
\left (\begin{smallmatrix}* 
& 0 & 0 & 0 & 
\ldots & 0 & ? & ? & ? & \ldots & ?\\ 0 & * & 0 & 0 & \ldots & 0 & ? & 
? & ? & \ldots & ?\\ 0 & 0 & * & 0 & 
\ldots & 0 & ? & ? & ? & \ldots & ?\\ 
\vdots &\vdots &\vdots & \vdots & \ddots & \vdots  & \vdots & \vdots & \vdots & \vdots & \vdots \\
0&0 &0 &0 &\ldots & * & ? &? &? &\ldots &? \\ 
\end{smallmatrix}\right )\bmod p$$
where each $*$ is coprime to $p$. 
\end{claim}

\noindent
The claim is an elementary fact from linear algebra. For convenience we give a proof at the end of 
the section. In the following let $p=2$ and let $\hat T= S^1\times \ldots \times S^1\to T$ denote a finite covering action 
of odd degree with the properties described in Claim \ref{first claim}. 

The next claim gives information on the action of involutions of $\hat T$ on $T_{pt} M$ (again the 
proof is postponed to the end of the section). This information will be used to exhibit involutions 
for which the fixed point manifold is connected of codimension $4$ (since $\hat T\to T$ has odd 
degree the same property holds for involutions of $T$). 

\begin{claim}\label{second claim}
\begin{enumerate}
\item For the involution $\sigma _i\in \hat T $ corresponding to $c_i$ 
($\sigma _i$ acts on $T_{pt}M$ by $\exp (c_i/2)\in U(1)^{2k}$) the component of $M^{\sigma _i}$ which 
contains the $\hat T$-fixed point $pt$ has codimension $4$, i.e. $c_i$ has two odd entries. 
\item The involution $\sum _i \sigma _i$ acts trivially on each representation space 
$R_i$, $i>2r$, i.e. each of the last $(2k-2r)$ columns of $A$ has an even number of odd entries.
\item The $2$-torus of $\hat T$ (i.e. the subgroup generated by the involutions $\sigma _i$) acts non-trivially on at most $r$ of the
 representation spaces $R_{2r+1},\ldots , R_{2k}$, i.e. at 
most $r$ of the last $(2k-2r)$ columns of $A$ are non-zero modulo $2$. 
\end{enumerate}
\end{claim}

\noindent
Hence, for every involution $\sigma \in \hat T$ and every $\hat T$-fixed point the connected 
component of $M^\sigma $ which contains the fixed point has codimension $\leq 6r\leq 
\frac 1 2 
\dim M$. Since $M$ has \rfpd \ for the prime $2$ it 
follows that $M^\sigma $ is the union of a connected component of codimension $\leq 6r$ with $\hat 
T$-fixed point (in fact the codimension is $\leq 4r$ by Theorem \ref{involution}) and (a possible empty set of) components with fixed point free $T$-action. Recall that $\codim M^{\sigma _i}=4$. Since any isolated 
$\sigma $-fixed point is also fixed by $\hat T$ it follows that $M^{\sigma _i}$ is connected of 
codimension $4$ for every $i$. 

Let $S_j$ denote the $j$th $S^1$-factor of $\hat T$. Below we shall use Theorem \ref{theorem cyclic} 
to show 

\begin{claim}\label{third claim} The fixed point manifold $M^{S_j}$ is 
connected of codimension $4$. In particular, $c_j$ has two odd entries and all other entries vanish. 
\end{claim}

\noindent
Note that Claim \ref{first claim}, Claim \ref{second claim} and Claim \ref{third claim} imply that the codimension of any connected component of $M^{\hat T}$ is $\leq 6r$.

With this information at hand we shall now complete the proof of Lemma \ref{mlemma}. By the second 
part of Claim \ref{second claim} one can choose $c_j$ such that 
$$\left (\begin{smallmatrix}c_{j}\\ c_{2r}\end{smallmatrix}\right )=
\left (\begin{smallmatrix}0 & \ldots & 0 & odd &0 & \ldots & 0 & 0 &0 & \ldots & 0 & odd & 0 & \ldots & 0\\
0 & \ldots & 0 &0 & 0 & \ldots & 0 & odd &0 & \ldots & 0 & odd & 0 & \ldots & 0\end{smallmatrix}\right )$$
Consider the $S^1$-subgroup $S\subset \hat T$ determined by $2\cdot c_j+c_{2r}\in \hat h (I_{\hat T})$. The fixed point manifold $M^\sigma $ of the involution $\sigma 
\in S$ is equal to $M^{\sigma _{2r}}$ which is connected of codimension $4$. The ${S}$-fixed point manifold has codimension $2$ in $M^\sigma $. Since $M^\sigma $ has \rfpd \ for the prime $2$ the fixed point manifold $M^S$ is connected. Also, $M^S$ is orientable. Since $M$ is a $Spin$-manifold $M^\sigma $ is orientable as well (see for example \cite{BoTa}, Lemma 10.1). Finally note that the same properties hold for the $S^1$-subgroup of $T$ which is the image of $S$ under the covering homomorphism $\hat T\to T$.\proofend 

\bigskip
\noindent
{\bf Proof of Claim \ref{first claim}:} It is an elementary fact from linear algebra that 
$h(I_T)\subset \Z 
^{2k}$ admits a basis $b_1,\ldots ,b_{2r}$ such that
$$\left (\begin{smallmatrix}b_1\\ 
\vdots \\ b_{2r}\end{smallmatrix}\right )=\left (\begin{smallmatrix}? & ? & ? & ? & 
\ldots & ? & ? & ? & ? & \ldots & ?\\ 0 & ? & ? & ? & \ldots & ? & ? & 
? & ? & \ldots & ?\\ 0 & 0 & ?& ? & 
\ldots & ? & ? & ? & ? & \ldots & ?\\ 
\vdots &\vdots &\vdots & \vdots & \ddots & \vdots  & \vdots & \vdots & \vdots & \vdots & \vdots \\
0&0 &0 &0 &\ldots & ?& ? &? &? &\ldots &? \\ 
\end{smallmatrix}\right ) 
\equiv 
\left (\begin{smallmatrix}* 
& 0 & 0 & 0 & 
\ldots & 0 & ? & ? & ? & \ldots & ?\\ 0 & * & 0 & 0 & \ldots & 0 & ? & 
? & ? & \ldots & ?\\ 0 & 0 & * & 0 & 
\ldots & 0 & ? & ? & ? & \ldots & ?\\ 
\vdots &\vdots &\vdots & \vdots & \ddots & \vdots  & \vdots & \vdots & \vdots & \vdots & \vdots \\
0&0 &0 &0 &\ldots & * & ? &? &? &\ldots &? \\ 
\end{smallmatrix}\right )\bmod p$$
after permuting columns (i.e. after permuting the representation spaces $R_i$) if necessary. Here 
each $*$ is coprime to $p$. Note that the choice of a basis of $h(I_T)$ is equivalent to the choice 
of an isomorphism $S^1\times \ldots \times S^1\to T$. Using suitable row operations of the form 
$b_i\rightsquigarrow \alpha \cdot b_i +\beta \cdot b_j$, where $i<j$, $\alpha \not \equiv 0 
\bmod p$ and $\beta \equiv 0 \bmod p$, one obtains a basis $c_1,\ldots , c_{2r}$ of $h(I_T)\otimes \Q 
$ such that the matrix $A=(c_1,\ldots ,c_{2r})^t$ has the properties given in the claim. Each $c_j$ determines a homomorphism $\iota _j:S^1\to T$. Since $c_1,\ldots ,c_{2r}$ is a basis  of $h(I_T)\otimes \Q $ the homomorphism $\hat T:=S^1\times \ldots \times S^1\overset {\iota _1\cdot \ldots \cdot \iota _{2r}}\longrightarrow T$ is a finite covering homomorphism. The induced homomorphism $\hat h: I_{\hat T}\to \Z ^{2k}$ of integral lattices maps the standard basis to $c_1,\ldots ,c_{2r}$. In view of the properties of the matrix $A$ the covering homomorphism has degree coprime to $p$.
\proofend

\bigskip
\noindent
{\bf Proof of Claim \ref{second claim}:} First note that by Theorem \ref{involution} and Corollary 
\ref{oddaction} for every involution $\sigma 
\in \hat T$ the fixed point manifold $M^\sigma $ has codimension $\leq 
4r$ and the dimension of each connected component of $M^\sigma $ is divisible by $4$. Since the 
dimension of $M$ is $\geq 8r$ and $M$ has \rfpd \ for the prime $2$ it follows that for every connected component $F$ 
of $M^\sigma $ either $\codim F\leq 4r$ or $\dim F\leq 4r-4$. 

Let $\hat h_1:I_{\hat T}\otimes (\Zp 2) \to (\Zp 2)^{2k}$ denote the mod $2$ reduction of $\hat h$. 
For the binary linear code ${\cal C}_1:=im(\hat h_1)\subset (\Zp 2)^{2k}$ we conclude that each code 
word $\sigma \in {\cal C}_1$ has either weight\footnote{The weight of $\sigma $ is defined as the 
number of entries equal to $1$.} $\wt (\sigma )\leq 2r$ or co-weight $\cowt (\sigma ):=2k-\wt (\sigma 
)\leq 2r-2$. In particular, the mod $2$ reduction of $c_i$, denoted by $\sigma _i$, has weight $\wt 
(\sigma 
_i)\leq 2r$. Since the weight function is sublinear, i.e. $\wt (\sigma 
+\sigma 
^\prime )\leq \wt (\sigma )+\wt (\sigma^\prime )$, and $2k\geq 6r$ it follows
that the subset of code words with weight $\leq 2r$ is closed under addition. Hence $\wt (\sigma 
)\leq 2r$ for every $\sigma \in {\cal C}_1$. In particular, this inequality holds for $\sum _j \sigma 
_j $ and $\sum 
_{j\neq i} \sigma _j$ which implies $\wt (\sigma _i)=2$ (i.e. $c_i$ has two odd entries)
 and implies that each of the last $(2k-2r)$ columns of $A$ has an even number of 
odd entries. Finally note that if $l$ of the last $(2k-2r)$ columns of $A$ are non-zero modulo $2$ then 
the weight of some code word is $\geq 2l$. Hence $l\leq r$. This completes the proof of the claim. 
\proofend 

\bigskip
\noindent
{\bf Proof of Claim \ref{third claim}:} To show that $M^{S_j}$ is connected of codimension $4$ it 
suffices to prove this property for all cyclic subgroups $C_s\subset S_j$ of order $2^s$, $s\in \N$. 
We know already that $M^{C_1}$ is connected of codimension $4$.

Assume $M^{C_s}$ is connected of codimension $4$. To show the corresponding property for $C_{s+1}$ we will use the following consequence of Theorem \ref{theorem cyclic} which we prove first: 

\bigskip
\noindent
{\it Let $M$ be a connected $Spin$-manifold with $S^1$-action. Assume the cyclic subgroup $C_s\subset 
S^1$ of order $2^s$ acts non-trivially on $M$ and the fixed point manifold $M^{C_s}$ is connected. If 
the first $(r+1)$ coefficients in the expansion $\Phi _0(M)$ do not vanish then the codimension of 
the submanifold $M^{C_{s+1}}$ of $M^{C_s}$ is $\leq 4r-2$.} 

\bigskip
\noindent
 Let $\oorder :=2^{s+1}$. Note that for some connected component $Y$ of the $S^1$-fixed 
point manifold $m_\oorder (Y)\leq r$ by Theorem \ref{theorem cyclic}. Let $Z\subset M^{C_{s+1}}$ be 
the connected component which contains $Y$. Note that $m_\oorder (Y)$ is strictly larger than $\codim 
(Z\subset M^{C_s})/4$ since $M^{C_s}$ is a proper connected submanifold of $M$. Hence the codimension 
of $Z$ in $M^{C_s}$ must be less than $4r$ which implies that the codimension of the submanifold 
$M^{C_{s+1}}$ of $M^{C_s}$ is $\leq 4r-2$.

\bigskip
\noindent
We shall now continue with the study of the action of $S_j$. By the statement above the codimension of the submanifold $M^{C_{s+1}}$ of $M^{C_{s}}$ is $\leq 
4r-2$. Since $M$ has \rfpd \ for the prime $2$ a connected component of $M^{C_{s+1}}$ has either codimension $\leq 
4r-2$ in $M^{C _{s}}$ or has dimension $\leq 4r-4$. 

Next consider the action of the $S^1$-subgroup $\widetilde S$ of $\hat T$ which is determined by $\widetilde 
c:=c_{j}+2^s\cdot \sum 
_i c_i\in \hat h(I_{\hat T})$. Let $\widetilde C_{s}\subset \widetilde S$ denote the cyclic subgroup of order $2^s$. By the statement above the codimension of the 
submanifold $M^{\widetilde C_{s+1}}$ of $M^{\widetilde C_{s}}$ is $\leq 4r-2$.

It follows that either the number of entries of $c_j$ which are $\equiv 2^s\bmod 2^{s+1}$ is 
$\leq 2r-1$ or the number of entries of $c_j$ which are $\equiv 0\bmod 2^{s+1}$ is $\leq 2r-2$. The 
same property holds for $\widetilde c$. Also the mod $2^{s+1}$-reductions of $c$ and $\widetilde c$ have the same last $(2k-2r)$ entries by the second part of Claim \ref{second claim}.

This implies that the mod $2^{s+1}$-reduction of $c_{j}$ has only two non-zero entries. In other 
words $C_{s+1}\subset S_j$ acts trivially on the tangent bundle of $M^{C_s}$ at the $T$-fixed point 
$pt$. Thus $M^{C_{s+1}}=M^{C_s}$ and $M^{C_{s+1}}$ is connected of codimension $4$. This completes 
the induction step. It follows that the fixed point manifold $M^{S_j}$ is connected of codimension 
$4$. 
\proofend

\section{Ricci versus sectional curvature}\label{versus}
Apparently the only known topological property which allows to distinguish positive Ricci curvature from 
positive curvature on simply connected manifolds is based on Gromov's Betti number theorem 
\cite{Gr2}: With respect to any field of coefficients the sum of Betti numbers of a positively curved 
$n$-dimensional manifold is less than a constant $C(n)$ depending only on the dimension. 

The bound in \cite{Gr2} which depends doubly exponentially on the dimension $n$ has been improved by 
Abresch \cite{Ab} who showed that the Betti number theorem holds for a bound $C(n)$ depending 
exponentially on $n^3$ (the bound $C(n)$ may be chosen to satisfy $\exp (5n^3+8n^2+4n+2)\leq C(n)\leq 
\exp (6n^3+9n^2+4n+4)$). Sha and Yang \cite{ShYa} constructed metrics of 
positive Ricci curvature on the $k$-fold connected sum of $S^n\times S^m$, $n,m>1$, for every $k\in 
\N $. By the Betti number theorem these manifolds do not admit positively curved metrics if $k$ is 
sufficiently large. 

In this section we present two methods to exhibit manifolds with ``small'' Betti numbers (i.e. the 
sum is less than $C(n)$) which admit metrics of positive Ricci curvature but do not admit metrics of 
positive curvature under assumptions on the symmetry. The first method which relies on Theorem 
\ref{main theorem} leads to 

\begin{theorem}\label{Ricci versus sec theorem1}
For every $r\geq 3$ and every $d\in \N $ there exists a simply connected closed manifold $M$ of 
dimension greater than $d$ such that: 
\begin{enumerate}  
\item $M$ admits a metric of positive Ricci curvature with symmetry rank $\geq 2r$.
\item $M$ does not admit a metric of positive curvature with symmetry rank $\geq 2r$.  
\end{enumerate}  
\end{theorem}  

\bigskip
\noindent
The manifold $M$ may be chosen to have small Betti numbers. The second method which relies on Theorem 
\ref{theorem cyclic} and recent work of Wilking (see Theorem \ref{connectivity theorem}) gives 
stronger information. Recall from the introduction that an $S^1$-action is called \almostisometric \ 
of order $\oorder $ if the cyclic subgroup of order $\oorder $ acts effectively and isometrically. 

\begin{theorem}\label{Ricci versus sec theorem2}
For every $d\in \N $ and every $\oorder \geq 2 $ there exists a simply connected closed manifold $M$ 
of dimension greater than $d$ such that: 
\begin{enumerate} 
\item $M$ admits a metric of positive Ricci curvature with \almostisometric \ $S^1$-action of order $\oorder$. 
\item $M$ does not admit a metric of positive curvature with \almostisometric \ $S^1$-action of order $\oorder$. 
\end{enumerate}
\end{theorem}

\bigskip
\noindent
Again the manifold $M$ may be chosen to have small Betti numbers. The examples occurring in both 
theorems are given by Riemannian products of the form $M_1\times M_2$. The first factor $M_1$ is a $Spin$-manifold of positive Ricci curvature with large symmetry rank and non-vanishing 
elliptic genus (below we shall choose $M_1$ to be a product of quaternionic projective planes). The second factor $M_2$ is a $Spin$-manifold of positive Ricci curvature for which the 
index $\hat A(M_2,TM_2)$ does not vanish. The next proposition shows that these properties hold for a 
hypersurface in $\C P^{n+1}$ of degree $n=2m>2$. 

\begin{proposition}\label{elliptic hypersurface proposition}
A non-singular hypersurface $V_{n}$ in $\C P^{n+1}$ of degree $n=2m>2$ has the following properties: 
\begin{enumerate}
\item $V_{n}$ is $Spin$ and admits a metric of positive Ricci curvature.
\item The Betti numbers of $V_n$ are small.
\item The index $\hat A(V_{n},TV_{n})$ does not vanish.
\end{enumerate}
\end{proposition}
\bigskip
\noindent
{\bf Proof:} For a proof of some of the properties of hypersurfaces used below we refer to \cite{Hi1} 
and references therein. Let $V^{l}$ denote a non-singular hypersurface in $\C P^{n+1}$ of degree $l$. 
The tangent bundle of $V^{l}$ is stably the restriction of $(T\C P^{n+1}-\gamma ^l)$ to $V^{l}$, 
where $\gamma $ denotes the complex line bundle over $\C P^{n+1}$ with first Chern class $h\in H^2(\C 
P^{n+1};\Z )$ dual to $\C P^n\hra \C P^{n+1}$ (in the following we shall denote the restriction of 
$h$ to $V^{l}$ also by $h$). Note that $\langle h^n, \lbrack V^{l}\rbrack \rangle=l$. The total Chern 
class of $V^{l}$ is given by
$$c(V^{l})=(1+h)^{n+2}\cdot (1+l\cdot h)^{-1}.$$
In particular, $c_1(V^{l})=(n+2-l)\cdot h$. Taking $l=n$ we conclude that $V_n$ is a $Spin$-manifold 
with positive first Chern class. By Yau's solution \cite{Yau1} of the Calabi conjecture $V_n$ admits 
a metric of positive Ricci curvature. 

The Euler characteristic of $V_n$ is equal to the Chern number $\langle c_n(V_n),\lbrack V_n\rbrack 
\rangle $ which can be computed via the formula for $c(M)$ given above. In turns out that
$$\chi (V_n)=\frac {(n-1)^{n+2} -1} n + (n+2).$$
By the Lefschetz theorem on hyperplane sections the odd Betti numbers of $V_n$ vanish for any field 
of coefficients (cf. \cite{Hi1}, \S 22). Hence the Betti numbers of $V_n$ are small. 

Finally we compute $\hat A(V_n,TV_n)$. By the cohomological version of the index theorem 
\cite{AtSiIII} $\hat A(V^{l}, TV^{l})$ is equal to $\langle B, \lbrack V^{l}\rbrack \rangle $, where 
$$B:=\left (\frac h {e^{h/2}-e^{-h/2}}\right )^{n+2}
\cdot \frac {e^{lh/2}-e^{-lh/2}} {lh} \cdot \left ( (n+2)e^h-e^{lh}\right ).$$
It follows that $\langle B, \lbrack V^{l}\rbrack \rangle $ is the residue of $l\cdot 
B/h^{n+1}$ at $h=0$. Changing variables, $w:=e^h-1$, one computes that $\hat A(V^{l}, TV^{l})$ is 
equal to the coefficient of $w^{n+1}$ in 
$$(1+w)^{(n-l)/2}\cdot ((1+w)^l-1)\cdot ((n+2)\cdot (1+w)-(1+w)^l).$$
Hence, $\hat A(V_{n}, TV_{n})=n+2-\binom {2n} {n+1}$ which does not vanish for $n=2m>2$. This 
completes the proof of the proposition.\proofend 

\bigskip
\noindent
{\bf Proof of Theorem \ref{Ricci versus sec theorem1}:} Let $M_1$ be the Riemannian product of 
$(r-1)$-copies of $\H P^2$, let $M_2$ be a non-singular hypersurface in $\C P^{n+1}$ of degree 
$n=2m\gg 0$ (equipped with a metric of positive Ricci curvature) and let $M$ be the Riemannian product 
$M_1\times M_2$. Since the symmetry rank of $\H P^2$ is three and $r\geq 3$ the symmetry rank of the 
positively Ricci curved $Spin$-manifold $M$ is $\geq 2r$. 

By Proposition \ref{elliptic hypersurface proposition} $\hat A(M_2,TM_2)\neq 0$. Thus the expansion $\Phi _0(M_2)$ has a pole of order $(\dim M_2)/8-1=(\dim M)/8-r$. 
For a homogeneous $Spin$-manifold the elliptic genus is equal to the signature \cite{HiSl}. In 
particular, $\Phi(\H P^2)=1$ which can also be shown by a direct computation. This gives
$$\Phi (M)=\Phi (\H P^2)^{r-1}\cdot \Phi (M_2)= \Phi (M_2).$$
It follows that $\Phi _0(M)$ has a pole of order $(\dim M)/8-r$. Thus, the first $(r+1)$ 
coefficient in the expansion $\Phi_0(M)$ do not vanish. By Theorem \ref{main theorem} the 
manifold $M$ does not admit a metric of positive curvature with symmetry rank $\geq 2r$.\proofend 

\bigskip
\noindent
{\bf Proof of Theorem \ref{Ricci versus sec theorem2}:} Let $M$ be the Riemannian product of $\H P^2$ 
and $V_n$ (the latter shall be equipped with a metric of positive Ricci curvature). Since the 
positively Ricci curved manifold $M$ has symmetry rank $\geq 3$ it admits a 
\almostisometric 
\ $S^1$-action of any order. 

Now assume $M$ carries a metric of positive curvature with \almostisometric \ $S^1$-action of order 
$\oorder $. We shall derive a contradiction for $n=2m\gg \oorder $. Let $\sigma \in S^1$ be of order 
$\oorder $. Since $\Phi _0(M)$ has a pole of order $(\dim M)/8-2$ the codimension of $M^\sigma $ is 
bounded from above by a constant which depends on $\oorder $ but does not depend on $\dim M=2n+8$ 
(see Corollary \ref{corollary special case}). 

Let $N\subset M^\sigma $ be a fixed point component of minimal 
codimension $k$. By Theorem \ref{connectivity theorem} the inclusion $N\hra M$ is 
$((2n+8)-(k+l))$-connected, where $l=k-1$. Using Poincar\'e duality for $M$ and $N$ it follows that 
$H^i(M;\Z )$ and $H^{i+k}(M;\Z )$ are isomorphic in the range $l<i<(2n+8)-(k+l)$. Hence, for 
$n\gg \oorder $ the cohomology groups $H^{n-k}(M;\Z )$ and $H^{n}(M;\Z )$ are isomorphic. 

Recall from \cite{Hi1}, \S 22, that the even Betti numbers $b_{2i}(V_n)$ are one except for $2i=n$ 
and the odd Betti numbers vanish. This implies $b_{n-k}(M)\leq 3$. However $b_{n}(M)$ is $>3$ 
(compare with the formula for $\chi (V_n)$ in the proof of the proposition above). This gives the 
desired contradiction. 
\proofend

\bigskip
\noindent
{\it Acknowledgements:} I like to thank Wilderich Tuschmann, Burkhard Wilking and Wolfgang Ziller for many 
stimulating discussions.

\bigskip
\noindent
Anand Dessai\\ e-mail: dessai@math.uni-augsburg.de\\ 
http://www.math.uni-augsburg.de/geo/dessai/homepage.html\\ Department of Mathematics, University of 
Augsburg, D-86159 Augsburg\\

\end{document}